# A general nonlinear mechanical model based on special-shaped cam-spring mechanism

Qiangqiang Li, Qingjie Cao

*Centre for Nonlinear Dynamics Research, School of Astronautics, Harbin Institute of Technology, Harbin 150001, China*

**Abstract**

Nonlinear systems exist widely in nature, however, how to construct systems with accurate expected non-linearity artificially is still a problem, which greatly limits their experimental study and engineering application. In this paper, we present a general nonlinear mechanical model (GNMM) based on special-shaped cam-spring mechanism (SCSM). GNMM consists of a mass, a roller, a track and a novel general spring model (GSM). Firstly, the equation of motion is derived and compared with that of SCSM. Then, the roller trajectory functions for arbitrary given continuous force are derived. Finally, Duffing system with softening non-linearity is constructed as a design example. The results presented herein show that GNMM has a more general restoring force than SCSM, and prove that there are at least six independent mechanical configurations for a given nonlinear system with continuous restoring force.

# 1. Introduction

Nonlinear systems exist widely in nature, however, how to construct them with accurate expected non-linearity artificially remains a problem, which limit the experimental study and engineer application of nonlinear systems. In recent years, owning to that the nonlinear characteristic directly determines the performance of the mechanical systems applied in the fields such as vibration isolation, vibration energy harvesting and robot engineering, designing mechanical system with expected non-linearity has received much attention[1-4].

Up to now, there are two main ways to design mechanical systems with expected non-linearity. For the first method, a linear spring is connected with a string wrapping around a non-circular pulley which rotates linearly, leading to a nonlinear torque[1]. By designing the pulley profile, we can obtain the expected nonlinear torque. This work can date back to Michel Jean, a French scientist, firstly uses this method to construct an accurate mechanical structure of Duffing system in 1965[5]. For the second method, a linear spring is connected with a roller which can slide along a special-shaped cam. With the cam moving in the direction perpendicular to the deformation direction of the spring, we can have a nonlinear restoring force, and it can be designed accurately by giving a corresponding shaped cam[2-4]. Collectively, the key constructing idea of these two methods is introducing a nonlinear boundary condition to linear spring to produce expected non-linearity. By contrast, the mechanism designed by the first method can only produce uni-directional torque since the string can only be pulled, using the second method to construct expected nonlinear mechanical system has become prevailing.

The motivations and contributions of this paper are (i) to present a more general nonlinear mechanical model (GNMM) by creatively introducing a general spring model (GSM) (with either positive or negative stiffness) to special-shaped cam-spring mechanism (SCSM), and (ii) to prove that there are at least six mechanical configurations for a nonlinear system with continuous restoring force.

The rest of the paper is organized as follows. In Sec.2, GSM is constructed and GNMM is constructed based on GSM. In Sec. 3, the inverse problem calculating corresponding roller trajectory functions for arbitrary given continuous force are proposed and solved. In Sec. 4, Duffing system with softening non-linearity is constructed as a design example. Finally, some conclusions of this study are drawn in Sec. 5.

# 2. General nonlinear mechanical model

## 1.1. GSM

Fig. 1 shows the model of GSM constructed by connecting a linear spring (positive stiffness element) and a dipteran flight mechanism[6] (a negative stiffness element) in parallel. GSM consists of a vertical spring with stiffness $K_1$, a pair of rigid rods with length $L$ and a pair of oblique springs with stiffness $K_2$ and original length $L_0$. The whole system is initially at the equilibrium position with vertical spring unstressed and the two rods with no angle about the horizontal direction, as shown in Fig. 2(a). With an external force applied to roller $M$ in $Y$ direction, as shown in Fig. 2(b), the roller deviates from the equilibrium position with a displacement $Y$. When $Y=L$, the roller will be locked by the two rods, as illustrated in Fig. 2(c), so the effective working range of GSM is limited within $|Y|<L$.

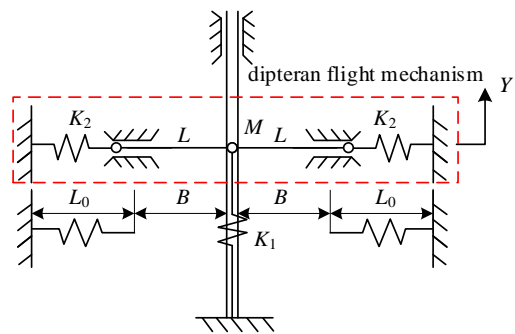

Fig. 1. The schematic diagrams of GSM.

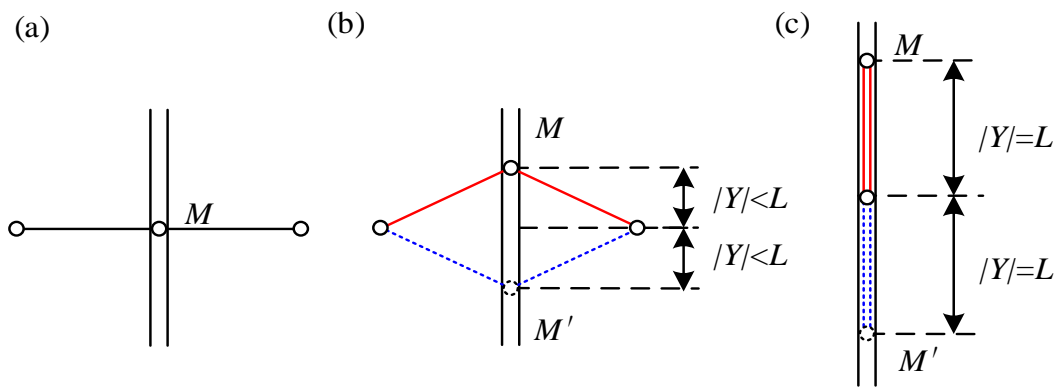

Fig. 2. Three different static states of GSM: (a) at equilibrium position; (b) deviation from equilibrium position with a displacement $Y$; (c) at 'locked' position.

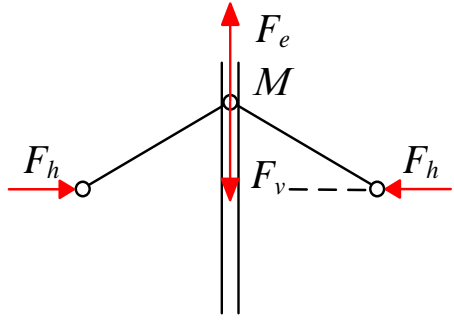

Fig. 3. Force analysis diagram of GSM

According to force analysis illustrated in Fig. (3), the relationship between the applied force and the displacement can be derived as

$$F_e = K_1 Y - 2K_2\left(1 - \frac{B}{\sqrt{L^2 - Y^2}}\right) Y, |Y| < L \tag{1}$$

where $B$ is half length between the inner ends of the horizontal springs at their equilibrium positions. By differentiating the first expression of Eq. (1) with respect to displacement $Y$, the stiffness of GSM is given as

$$K_{GSM} = K_1 - 2K_2 + \frac{2K_h B L^2}{(L^2 - Y^2)^{3/2}}, |Y| < L. \tag{2}$$

When $B = 0$, Eq. (3) can be rewritten as

$$K_{GSM} = K_1 - 2K_2. \tag{3}$$

It can be seen from Eq. (3) that the stiffness of GSM is determined by two parameters: the vertical spring stiffness $K_v$ and the oblique spring stiffness $K_h$. By adjusting the two parameters we can have negative, zero or positive linear stiffness, which means GSM can work as a general linear spring. Compared with the traditional spring, the difference lies in when GSM works as a negative spring, the deformation direction is the same as that of the traditional spring, but the spring force acts in the opposite direction.

**1.2. GNMM**

Based on the general spring GSM, GNMM is established, as shown in Fig. 4(b), which comprises a mass, a track (consolidated with the mass) and GSM. Compared with SCSM as shown in Fig. 4(a), the basic geometrical configuration (placing a nonlinear boundary to a linear spring) has not changed, the improvements are: 1) the cam is abstracted to a track that can not only be compressed but also be stressed and 2) the traditional linear spring is substituted by GSM with either positive or negative stiffness.

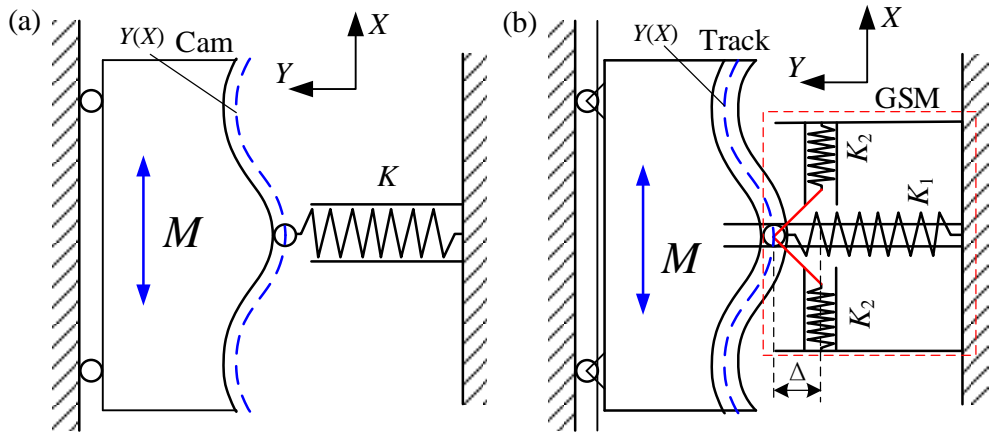

Fig. 4. Physical models of (a) SCSM given by [3] and (b) GNMM

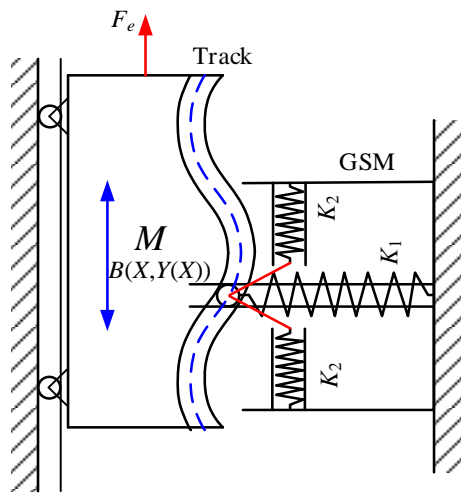

Fig. 5. GNMM deviates form the initial position with an external force.

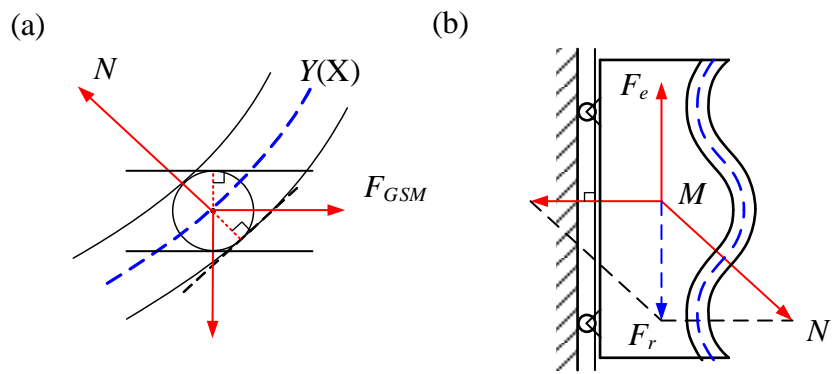

Fig. 6. Force analysis diagrams of (a) the roller and (b) the mass.

To give a clear description of the static characteristic of GNMM, an external force $F_e$ is applied to the lumped mass in $X$ direction. Without the applied force, as shown in Fig. 4(b), GNMM stays at the initial position with GSM pre-stretched (pre-compressed state is also permitted) by length Δ. Under the applied force, the mass moves by a displacement $X$, and the roller roles from the initial point $A$ to point $B$, which stretches (or compresses) GSM by a displacement $Y(X)$, as illustrated in Fig. 5. In this process, though GSM works as a linear spring, a nonlinear spring force $F_{GSM}$ with respect to displacement $X$ happens in $Y$ direction due to the nonlinear track, as illustrated by Fig. 6(a). The nonlinear track plays a vital role not only in no-linearizing GSM but also in converting the nonlinear spring force to the motion direction by the reacting force $N'$ of supporting force $N$, as illustrated Fig. 6(b), which results in a nonlinear restoring force $F_r$ equalling to the applied force but in opposite direction.

Based on the static analysis of GNMM, under the applied force, the nonlinear spring force of GSM with respect displacement $X$ can be given as

$$F_{GSM}(X) = -K_{GSM} Y(X) \tag{4}$$

By converting $F_{GSM}$ to the motion direction, the restoring force of GNMM can be derived and written as

$$F(X) = -K_{GSM} Y(X) \frac{dY(X)}{dX} \tag{5}$$

Assuming the damping and friction of the system are not considered, the equation of motion of GNMM with no external excitation is derived as

$$M\ddot{X} + K_{GSM} Y(X) \frac{dY(X)}{dX} = 0 \tag{6}$$

with the definition domain written as

$$X \in \{X \mid -L < Y(X) < L\} \tag{7}$$

From Eq. (6), it can be seen that the restoring force of GNMM is a first-order differential expression about $Y(X)$ that must meet the existence of the first derivative in the definition domain. Compared with the restoring force of SCSM given by [6], that of GNMM is a more general form lies in that the stiffness parameter can be taken any real number (without limitation to be positive) and $Y(X)$ has symmetric definition range about $Y(X)=0$ (without limitation to [-Δ$_{max}$, 0] where Δ$_{max}$ is the maximum compression of traditional linear spring).

## 3. Roller trajectory function

In contrast to the constructing process of GNMM discussed above, the practical design process is an inverse problem solving an initial value problem of differential equation about $Y(X)$:

$$\begin{cases} -K_{GSM}Y(X)\dfrac{dY(X)}{dX} = F(X), \\ Y(0) = \Delta, \\ -L < Y(X) < L. \end{cases} \quad (8)$$

For a given continuous restoring force $F(X)$, the general solution can be given as.

$$Y_{1,2}(X) = \pm\sqrt{\Delta^2 - \dfrac{2}{K_{GSM}}\int_0^X F(X)dX}, X \in \left\{ X \mid \Delta^2 - L^2 < \dfrac{2}{K_{GSM}}\int_0^X F(X)dX < \Delta^2 \right\}. \quad (9)$$

In the case when $\Delta \neq 0$, according to the sign of $K_{GSM}$, $Y_{1,2}(X)$ can be divided into four different forms:

$$\begin{cases} Y_{11}(X) = \sqrt{\Delta^2 - \dfrac{2}{K_{GSM}}\int_0^X F(X)dX}, X \in \left\{ X \mid \Delta^2 - L^2 < \dfrac{2}{K_{GSM}}\int_0^X F(X)dX < \Delta^2, K_{RSS} > 0 \right\}. \\ Y_{12}(X) = \sqrt{\Delta^2 - \dfrac{2}{K_{GSM}}\int_0^X F(X)dX}, X \in \left\{ X \mid \Delta^2 - L^2 < \dfrac{2}{K_{GSM}}\int_0^X F(X)dX < \Delta^2, K_{RSS} < 0 \right\} \\ Y_{21}(X) = -\sqrt{\Delta^2 - \dfrac{2}{K_{GSM}}\int_0^X F(X)dX}, X \in \left\{ X \mid \Delta^2 - L^2 < \dfrac{2}{K_{GSM}}\int_0^X F(X)dX < \Delta^2, K_{RSS} > 0 \right\} \\ Y_{22}(X) = -\sqrt{\Delta^2 - \dfrac{2}{K_{GSM}}\int_0^X F(X)dX}, X \in \left\{ X \mid \Delta^2 - L^2 < \dfrac{2}{K_{GSM}}\int_0^X F(X)dX < \Delta^2, K_{RSS} < 0 \right\} \end{cases} \quad (10)$$

In the case when $\Delta = 0$, according to the sign of $K_{GSM}$, $Y_{1,2}(X)$ can be divided into four new different forms:

$$\begin{cases} Y_{13}(X) = \sqrt{-\dfrac{2}{K_{GSM}}\int_0^X F(X)dX}, X \in \left\{ X \mid -L^2 < \dfrac{2}{K_{GSM}}\int_0^X F(X)dX < 0, K_{RSS} > 0 \right\}. \\ Y_{14}(X) = \sqrt{-\dfrac{2}{K_{GSM}}\int_0^X F(X)dX}, X \in \left\{ X \mid -L^2 < \dfrac{2}{K_{GSM}}\int_0^X F(X)dX < 0, K_{RSS} < 0 \right\} \\ Y_{23}(X) = -\sqrt{-\dfrac{2}{K_{GSM}}\int_0^X F(X)dX}, X \in \left\{ X \mid -L^2 < \dfrac{2}{K_{GSM}}\int_0^X F(X)dX < 0, K_{RSS} > 0 \right\} \\ Y_{24}(X) = -\sqrt{-\dfrac{2}{K_{GSM}}\int_0^X F(X)dX}, X \in \left\{ X \mid -L^2 < \dfrac{2}{K_{GSM}}\int_0^X F(X)dX < 0, K_{RSS} < 0 \right\} \end{cases} \quad (11)$$

For solutions given by Eq. (10), it can be easily proved that they always exist for arbitrary continuous $F(X)$ by intermediate value theorem of continuous function. For solutions given by Eq. (11), under the condition that $\int_0^X F(X)dX < 0$ ($\int_0^X F(X)dX > 0$), $Y_{13}$ and $Y_{14}$ ($Y_{23}$ and $Y_{24}$) do not exist, and under other conditions, $Y_{13}$, $Y_{14}$, $Y_{23}$ and $Y_{24}$ exist simultaneously. Therefore, there are at least six and at most eight solutions to Eq. (8), which means that at least six mechanical configurations can be found in physical world based on GNMM for a given nonlinear system with continuous restoring force.

## 4. A design example: $M\ddot{X} - 5 \times 10^3 X^3 = 0$

The general expression of Duffing system with softening non-linearity can be written as

$$M\ddot{X} + AX^3 = 0 \tag{12}$$

Where $A$ is negative constant. Combining the restoring force of system (12) with Eq. (9) yields

$$Y_{1,2}(X) = \pm\sqrt{\Delta^2 + \frac{A}{2K_{GSM}}X^4}, X \in \left\{X \mid -\Delta^2 < \frac{A}{2K_{GSM}}X^4 < L^2 - \Delta^2\right\}. \tag{13}$$

which can be divided into six valid solutions:

$$\begin{cases} Y_{11}(X) = \sqrt{\Delta^2 + \frac{A}{2K_{GSM}}X^4}, X \in \left\{X \mid -\left(-\frac{2K_{GSM}\Delta^2}{A}\right)^{\frac{1}{4}} < X < \left(-\frac{2K_{GSM}\Delta^2}{A}\right)^{\frac{1}{4}}, K_{GSM} > 0\right\} \\ Y_{12}(X) = \sqrt{\Delta^2 + \frac{A}{2K_{GSM}}X^4}, X \in \left\{X \mid -\left(\frac{2K_{GSM}(L^2 - \Delta^2)}{A}\right)^{\frac{1}{4}} < X < \left(\frac{2K_{GSM}(L^2 - \Delta^2)}{A}\right)^{\frac{1}{4}}, K_{GSM} < 0\right\} \\ Y_{13}(X) = \sqrt{\frac{A}{2K_{GSM}}}X^2, X \in \left\{X \mid -\left(-\frac{2K_{GSM}L^2}{A}\right)^{\frac{1}{4}} < X < \left(-\frac{2K_{GSM}L^2}{A}\right)^{\frac{1}{4}}, K_{GSM} < 0\right\} \\ Y_{21}(X) = -\sqrt{\Delta^2 + \frac{A}{2K_{GSM}}X^4}, X \in \left\{X \mid -\left(\frac{2K_{GSM}\Delta^2}{A}\right)^{\frac{1}{4}} < X < \left(\frac{2K_{GSM}\Delta^2}{A}\right)^{\frac{1}{4}}, K_{GSM} > 0\right\} \\ Y_{22}(X) = -\sqrt{\Delta^2 + \frac{A}{2K_{GSM}}X^4}, X \in \left\{X \mid -\left(\frac{2K_{GSM}(L^2 - \Delta^2)}{A}\right)^{\frac{1}{4}} < X < \left(\frac{2K_{GSM}(L^2 - \Delta^2)}{A}\right)^{\frac{1}{4}}, K_{GSM} < 0\right\} \\ Y_{23}(X) = -\sqrt{\frac{A}{2K_{GSM}}}X^2, X \in \left\{X \mid -\left(-\frac{2K_{GSM}L^2}{A}\right)^{\frac{1}{4}} < X < \left(-\frac{2K_{GSM}L^2}{A}\right)^{\frac{1}{4}}, K_{GSM} < 0\right\} \end{cases}$$

$$\tag{14}$$

It is assumed that the nonlinear coefficient $A$ is set as -5000 Nm$^{-4}$, and the physical parameters $K_{GSM}$, $\Delta$ are selected as 100 Nm$^{-1}$/-100Nm$^{-1}$ and 0.1 m/-0.1 m, respectively. According to the roller trajectory functions given by Eq. (14), the corresponding mechanical configurations of GNMM are plotted in Fig. 7 (a-f), respectively. It can be seen that Model (a) is a special case of GNMM with $K_{GSM}>0$, $\Delta>0$, which corresponds to the mechanical configuration of SCSM. Benefiting from the generality stiffness of GSM, model (b) ($K_{GSM}<0$, $\Delta>0$) and model (c) ($K_{GSM}<0$, $\Delta=0$) can be constructed based on GNMM. In addition, owning to the potential energy of GNMM is an even function of $Y(X)$ (the same potential energy function means the same restoring force), model (e-f) which have symmetric roller trajectory with that of model (e-f) about $X$ axis, respectively can also be obtained. These six models are just like sextuplets in the physical word mirrored with the mathematical model of Duffing system: $M\ddot{X} - 5 \times 10^3 X^3 = 0$.

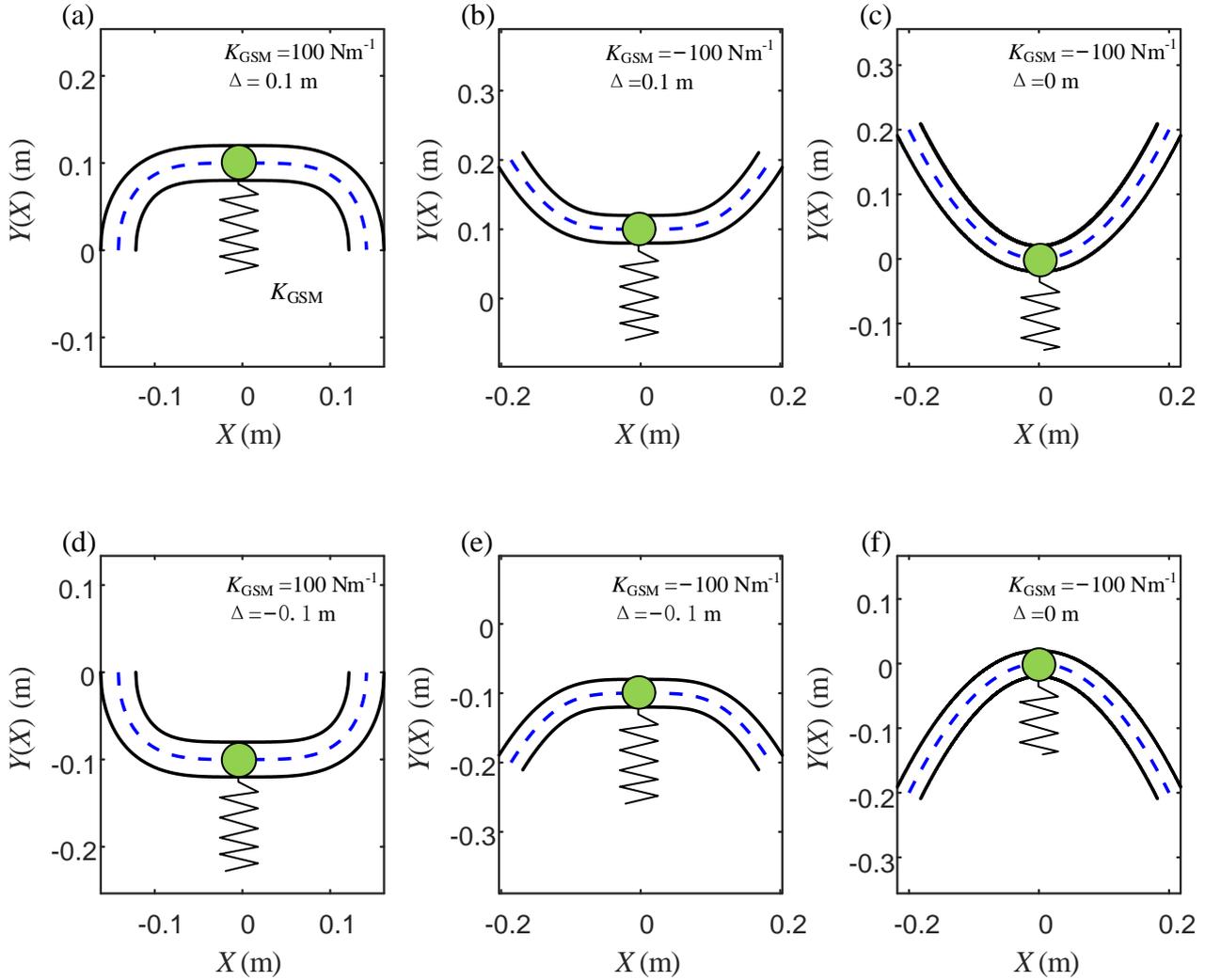

Fig. 7. Six different mechanical configurations of Duffing system: $M\ddot{X} - 5 \times 10^3 X^3 = 0$.

## 5. Conclusions

In this study, a general nonlinear mechanical model (GNMM) has been proposed by introducing GSM to SCSM. The equation of motion is derived and compared with that of SCSM. The roller trajectory functions for arbitrary given continuous force are derived. Duffing system with softening non-linearity is constructed as a design example. The main conclusions are drawn as:

(i) Compared with SCSM, GNMM is a more general model with the restoring force given by $F(X) = -K_{GSM} Y(X) \dfrac{dY(X)}{dX}, X \in \{X \mid -L < Y(X) < L\}$.

(ii) There are at least six mechanical configurations can be obtained based on GNMM for a given nonlinear system with continuous restoring force.

**Declaration of interests**

The authors declare that they have no known competing financial interests or personal relationships that could have appeared to influence the work reported in this paper.

**Author statement**

**Qiangqiang Li**: Methodology, Writing-original draft, Investigation

**Qingjie Cao**: Writing review and editing, Project administration, Resources

**Acknowledgments**

The authors would like thank professor Michel Jean for selflessly providing the English translation manuscript of his doctoral dissertation written in 1965.

**References**

[1] N. Schmit, M. Okada, Design and Realization of a Non-Circular Cable Spool to Synthesize a Nonlinear Rotational Spring, Advanced Robotics 26 (2012) 235-252. 10.1163/156855311x614545

[2] D. Zou, G. Liu, Z. Rao, T. Tan, W. Zhang, W.-H. Liao, A device capable of customizing nonlinear forces for vibration energy harvesting, vibration isolation, and nonlinear energy sink, Mechanical Systems and Signal Processing 147 (2021). 10.1016/j.ymssp.2020.107101

[3] D. Zou, G. Liu, Z. Rao, Y. Zi, W.-H. Liao, Design of a broadband piezoelectric energy harvester with piecewise nonlinearity, Smart Materials and Structures 30 (2021). 10.1088/1361-665X/ac112c

[4] S. Zuo, D. Wang, Y. Zhang, Q. Luo, Design and testing of a parabolic cam-roller quasi-zero-stiffness vibration isolator, International Journal of Mechanical Sciences 220 (2022). 10.1016/j.ijmecsci.2022.107146

[5] M. Jean, Sur les solutions périodique des équations différentielles de la mécaniques(PH.D. thesis), in, Aix-Marseille Université, 1965.

[6] Q. Cao, Y. Xiong, M. Wiercigroch, A novel model of dipteran flight mechanism, International Journal of Dynamics & Control 1 (2013) 1-11.